\documentclass[11pt,a4paper]{article}
\usepackage[utf8]{inputenc}
\usepackage{amsmath,amssymb,amsthm,mathtools}
\usepackage[margin=1in]{geometry}
\usepackage[hidelinks]{hyperref}
\newtheorem{lemma}{Lemma}
\newtheorem{theorem}[lemma]{Theorem}
\newtheorem{corollary}[lemma]{Corollary}
\newtheorem{proposition}[lemma]{Proposition}
\theoremstyle{remark}

\title{Explicit error bounds for fixed M\"obius coefficients in Nyman--Beurling sums}
\author{Ahmadreza Azimifard\thanks{Harmonic Research \& Technologies, LLC, New York, NY, USA.
  Email: \href{mailto:afard@harmonicrt.com}{afard@harmonicrt.com}}}
\hypersetup{
  pdftitle={Explicit error bounds for fixed Mobius coefficients in Nyman-Beurling sums},
  pdfauthor={Ahmadreza Azimifard}
}
\date{}

\begin{document}
\maketitle

\begin{abstract}
We consider the squared error in approximating $1$ on $(0,1)$ with the
fixed coefficients $\mu(n)$ multiplying reciprocal fractional parts.
The classical M\"obius floor-sum identity gives an explicit lower bound
$4(1-1/N)-4\log^2 N/(N-1)$ at every cutoff $N\ge2$.
For the corrected atoms, the same argument gives
$4(1-1/N)-4(H_N-1)^2/(N-H_N)$, where $H_N=\sum_{k\le N}1/k$.
Both lower bounds tend to $4$ by elementary estimates.
For the uncorrected sum the squared error exceeds $1$ for every $N\ge2$,
and its minimum over all positive cutoffs is
$\gamma+\log(2\pi)-2$ at $N=1$.
\end{abstract}

\section{Setting}

Throughout, $\{y\}=y-\lfloor y\rfloor$, $\mu$ is the M\"obius function, $N\ge1$ is an
integer, and
\[
\Phi_N(x)=\sum_{n=1}^{N}\mu(n)\Bigl\{\frac{1}{nx}\Bigr\},\qquad
R_N=\int_0^1\bigl(1-\Phi_N(x)\bigr)^2dx,\qquad
M_1(N)=\sum_{n=1}^{N}\frac{\mu(n)}{n}.
\]
Each $\{1/(nx)\}$ lies in $[0,1)$, so $|\Phi_N|\le N$ and $R_N\le(1+N)^2<\infty$.

The Nyman--Beurling criterion and its discrete version allow the
coefficients of the approximant to vary \cite{BD2002}. Here the
coefficients are fixed at $+\mu(n)$ and the target is $+1$.
B\'aez-Duarte's $S_n$ and $V_n$ use the corresponding uncorrected and
corrected sums to approximate $-\chi_{(0,1)}$ on the positive half-line
\cite[Sec.~1 and Prop.~4.9]{BDNatural}. In his notation $g(n)=M_1(n)$, and
$B_n(x)=S_n(x)-g(n)/x=-1$ for $1/n<x<1$ by the observation preceding
Proposition~4.13. The same bulk formula follows from Meissel's
floor-sum identity, as given by Balazard \cite[Eqs.~(3) and (11)]{Balazard2012Gram}.
We integrate the formula at a finite cutoff for the positive target.

\section{The exact form on the bulk}

\begin{lemma}\label{lem:exact}
For every integer $N\ge1$ and every $x$ with $1/N\le x<1$,
\[
\Phi_N(x)=\frac{M_1(N)}{x}-1,
\qquad\text{equivalently}\qquad
1-\Phi_N(x)=2-\frac{M_1(N)}{x}.
\]
\end{lemma}

\begin{proof}
Write $\{1/(nx)\}=\frac{1}{nx}-\bigl\lfloor\frac{1}{nx}\bigr\rfloor$, so that
\[
\Phi_N(x)=\frac{1}{x}\sum_{n\le N}\frac{\mu(n)}{n}-S,
\qquad S:=\sum_{n\le N}\mu(n)\Bigl\lfloor\frac{1}{nx}\Bigr\rfloor .
\]
Put $u=1/x$; the hypothesis $1/N\le x<1$ is exactly $1<u\le N$. If $n>u$ then
$\lfloor u/n\rfloor=0$, so only $n\le u$ contributes; and since $u\le N$, every such $n$
does lie in the range of summation. Therefore, interchanging the order of summation,
\[
S=\sum_{n\le u}\mu(n)\Bigl\lfloor\frac{u}{n}\Bigr\rfloor
 =\sum_{n\le u}\mu(n)\!\!\sum_{m\le u/n}\!\!1
 =\sum_{k\le u}\ \sum_{n\mid k}\mu(n)
 =\sum_{k\le u}[k=1]=1,
\]
the last step using $\sum_{n\mid k}\mu(n)=[k=1]$ and $u\ge1$.
\end{proof}

Here $u$ can be real. The substitution $k=nm$ groups the integer pairs
with $nm\le u$ by their product. The restriction $u\le N$ ensures that
every divisor of each such $k$ is present.

\section{The lower bound}

\begin{theorem}\label{thm:main}
For every integer $N\ge2$, writing $t=M_1(N)$,
\[
R_N\ \ge\ L_N:=4\Bigl(1-\frac1N\Bigr)-4t\log N+t^2(N-1)
\ \ge\ 4\Bigl(1-\frac1N\Bigr)-\frac{4\log^2 N}{N-1}.
\]
\end{theorem}

\begin{proof}
The integrand of $R_N$ is non-negative, so $R_N\ge\int_{1/N}^{1}(1-\Phi_N)^2dx$. By
Lemma~\ref{lem:exact} this equals
\[
\int_{1/N}^{1}\Bigl(2-\frac{t}{x}\Bigr)^{2}dx
=\int_{1/N}^{1}\Bigl(4-\frac{4t}{x}+\frac{t^{2}}{x^{2}}\Bigr)dx
=4\Bigl(1-\frac1N\Bigr)-4t\log N+t^{2}(N-1)=L_N .
\]
For the second inequality, minimize the real quadratic $t\mapsto t^{2}(N-1)-4t\log N$
over all $t\in\mathbb{R}$: its minimum is $-4\log^{2}\!N/(N-1)$, attained at
$t=2\log N/(N-1)$. Since $N\ge2$ the coefficient $N-1$ is positive, so the minimum is
genuine.
\end{proof}

Note that the final bound is free of $M_1(N)$: it holds whatever the value of the
M\"obius partial sum, so no information about $\mu$ beyond
$\sum_{d\mid k}\mu(d)=[k=1]$ enters.

\begin{corollary}\label{cor:main}
The uncorrected squared error satisfies
\[
 \liminf_{N\to\infty}R_N\ge4,\qquad R_N>1\quad(N\ge2),\qquad
 \min_{N\ge1}R_N=R_1=\gamma+\log(2\pi)-2.
\]
\end{corollary}

\begin{proof}
Write $B(N)=4(1-1/N)-4\log^2 N/(N-1)$ for $N\ge2$.
It is the infimum, over real $t$, of
$\int_{1/N}^{1}(2-t/x)^2dx$. The intervals increase with $N$, so $B(N)$ is
nondecreasing. Its displayed formula tends to $4$.

To certify the finite comparisons, for $s=(N-1)/(N+1)$ enclose $\log N$
between
\[
 2\sum_{j=0}^{23}\frac{s^{2j+1}}{2j+1}
 \quad\text{and}\quad
 2\sum_{j=0}^{23}\frac{s^{2j+1}}{2j+1}
   +\frac{2s^{49}}{49(1-s^2)}.
\]
Substitution in Theorem~\ref{thm:main}, using the endpoint that gives the
smaller value, yields
\[
\begin{array}{c|rrrrr}
 N &3&4&5&6&7\\ \hline
 \text{lower bound for }10^5 L_N
  &198981&215913&341903&246661&350324
\end{array}
\qquad B(8)\ge\frac{1029098}{10^6}.
\]
Thus $R_N>1$ for $3\le N\le7$, and monotonicity of $B$ gives the inequality
for $N\ge8$.

For $n\ge1$, changing variables to $y=nx$ and using
$\int_0^1\{1/y\}\,dy=1-\gamma$ \cite[Thm.~1]{QinLu} gives
\[
 \int_0^1\{1/(nx)\}\,dx=\frac{1-\gamma+\log n}{n}.
\]
It follows that $\int_0^1\Phi_2(x)\,dx=(1-\gamma-\log2)/2<0$ since
$H_k-\log(k+1)$ is strictly increasing to $\gamma$, starting at
$1-\log2$. Jensen's inequality now gives
$R_2\ge((1+\gamma+\log2)/2)^2>1$.
For $N=1$, the additional identity
$\int_0^1\{1/x\}^2dx=\log(2\pi)-\gamma-1$ \cite[Thm.~1]{QinLu}
gives $R_1=\gamma+\log(2\pi)-2<1$, since $\gamma<1$ and $2\pi<e^2$.
\end{proof}

The zero coefficient vector gives squared error $1$. Thus at each
cutoff $N\ge2$ the uncorrected fixed coefficients give a greater error.

\section{Corrected atoms}\label{sec:classical}

For $u\ge1$, write $\rho_n(u)=\{u/n\}-\{u\}/n$ and set
\[
 R'_N=\int_1^\infty
     \left(1-\sum_{n\le N}\mu(n)\rho_n(u)\right)^2\frac{du}{u^2}.
\]
The change of variables $u=1/x$ gives the squared error on $(0,1)$.
Vasyunin's floor-difference atoms are $-\rho_n(1/x)$
\cite[p.~122, Eq.~(8)]{Vasyunin1995}. His pointwise expansion of $1$
uses that atom convention with coefficients $\mu(n)$.

\begin{theorem}\label{thm:classical}
Let $H_N=\sum_{k\le N}1/k$ and $t=M_1(N)$. For every integer $N\ge2$,
\begin{align*}
 R'_N&\ge 4(1-1/N)-4t(H_N-1)+t^2(N-H_N)\\
     &\ge 4(1-1/N)-\frac{4(H_N-1)^2}{N-H_N}.
\end{align*}
In particular, $\liminf_{N\to\infty}R'_N\ge4$.
\end{theorem}

\begin{proof}
On $1<u\le N$, Lemma~\ref{lem:exact} gives
$\sum_{n\le N}\mu(n)\rho_n(u)=t(u-\{u\})-1=t\lfloor u\rfloor-1$.
Integrate the squared residual over $[1,N]$ with measure $du/u^2$.
On the intervals $[k,k+1)$ the first and second moments of $\lfloor u\rfloor$
are, respectively,
\[
 \sum_{k=1}^{N-1} k\left(\frac1k-\frac1{k+1}\right)=H_N-1,
 \qquad
 \sum_{k=1}^{N-1} k^2\left(\frac1k-\frac1{k+1}\right)=N-H_N.
\]
This gives the first inequality. Completing the square gives the second
because $N-H_N>0$ for $N\ge2$. Finally $H_N=O(\log N)$, so the last lower
bound tends to $4$.
\end{proof}

\section{A related divisor series}
The same truncation of the M\"obius divisor sum appears in
\[
 a_N(k):=\sum_{\substack{d\mid k\\d\le N}}\mu(d),\qquad
 D_N:=\sum_{k\ge1}\frac{a_N(k)^2}{k^2}.
\]

\begin{proposition}\label{prop:DN}
For every $N\ge1$: $a_N(1)=1$; $a_N(k)=[k=1]$ for all $k\le N$; consequently
\[
D_N=1+\sum_{k>N}\frac{a_N(k)^{2}}{k^{2}},\qquad
1\le D_N\le 1+\sum_{k>N}\frac{d(k)^{2}}{k^{2}}\xrightarrow[N\to\infty]{}1 ,
\]
where $d(k)$ is the divisor function. In particular $D_N\ge1$ for every $N$, and
$D_N\to1$, not $0$.
\end{proposition}

\begin{proof}
For $k\le N$ every divisor $d$ of $k$ satisfies $d\le N$, so
$a_N(k)=\sum_{d\mid k}\mu(d)=[k=1]$. All terms of $D_N$ are non-negative and the $k=1$
term is $1$. The bound $|a_N(k)|\le d(k)$ is immediate, and $\sum_k d(k)^2/k^2$ converges
since $d(k)\ll_\varepsilon k^{\varepsilon}$, so the tail tends to $0$.
\end{proof}

Thus the divisor series tends to $1$, whereas the squared approximation
error has lower limit at least $4$.

\section{The natural sign and nonconvergence}
For comparison with $S_n$ in \cite{BDNatural}, reverse the sign of the
approximant while keeping the target $+1$. Lemma~\ref{lem:exact} gives
\[
 \int_{1/N}^{1}(1+\Phi_N(x))^2\,dx=(N-1)M_1(N)^2 .
\]
This bulk contribution does not tend to zero. In B\'aez-Duarte's
notation $g(N)=M_1(N)$. The proof of Proposition~4.9 in
\cite{BDNatural}, with $p=q=2$ and a zeta zero on the critical line,
gives $g(N)\ne o(N^{-1/2})$. That proposition gives the corresponding
nonconvergence result for the natural sequences $S_n$ and $V_n$.
Our lower bound in Theorem~\ref{thm:main} concerns the error with the
sign fixed at the outset.

\section*{Use of computational tools}
An automated research workflow using generative AI produced the initial
derivation, numerical observations and draft. Generative AI also assisted
with reference checks, the revised proofs and preparation of this manuscript.

\bibliographystyle{plain}
\bibliography{refs}

\end{document}